\documentstyle[psfig,11pt]{pspum-l}

\title[Homology of $PGL_2$]{Integral Homology of $PGL_2$ over Elliptic Curves}
\author{Kevin P. Knudson}\thanks{Supported by an NSF Postdoctoral
Fellowship, grant no. DMS--9627503} 
\address{Department of Mathematics, Northwestern University, Evanston,
IL  60208}
\email{knudson@@math.nwu.edu}
\date{October 27, 1997}

\newtheorem{theorem}{Theorem}[section]
\newtheorem{prop}[theorem]{Proposition}

\newtheorem{cor}[theorem]{Corollary}

\newcommand{\X}{{\mathcal X}}
\newcommand{\D}{{\mathcal D}}
\newcommand{\zz}{{\Bbb Z}}
\newcommand{\cH}{{\mathcal H}}
\newcommand{\uk}{k^\times}
\newcommand{\ukw}{k(\omega)^\times}
\newcommand{\pgla}{PGL_2(A)}
\newcommand{\pglk}{PGL_2(k)}

\begin{document}
\maketitle

The Friedlander--Milnor Conjecture \cite{fried} asserts that if $G$ is
a reductive algebraic group over an algebraically closed field $k$,
then the comparison map 
$$H_{\text{\'et}}^\bullet(BG_k,\zz/p) \longrightarrow H^\bullet(BG,\zz/p)$$
is an isomorphism for all primes $p$ not equal to the characteristic
of $k$.  Gabber's rigidity theorem \cite{gabber} implies that this map
is indeed an isomorphism for the stable general linear group $GL$
(this is due to Suslin \cite{suslinloc} for $k={\Bbb C}$ and to
Jardine \cite{jardine} for arbitrary $k$).  Similarly, a proof of an
unstable version of rigidity would lead to a proof of the unstable
Friedlander--Milnor Conjecture.

In this note we consider unstable rigidity for the group $PGL_2$ over an
elliptic curve $E$.  We assume that $E$ is defined by the equation
$F(x,y)=0$, where $$F(x,y)=y^2 + a_1xy +a_3y - x^3 - a_2x^2 - a_4x
-a_6,$$ and the $a_i$ lie in an infinite field $k$. 
  Denote by $\overline{E}$ the projective curve $E \cup
\{\infty\}$. Denote by $A$ the coordinate ring of the affine curve
$E$.   If $l \in k$ and
$F(l,y)=0$ has no rational solutions, denote by $k(\omega)$ the
quadratic extension of $k$ inside the algebraic closure $\overline{k}$
for which $F(l,\omega)=0$.  Our main result is the following.

\medskip

{\sc Theorem.} {\em For all $i\ge 1$,
\begin{eqnarray*}
H_i(\pgla,\zz) & = & {\bigoplus\begin{Sb}
                          {p \in \overline{E}}\\ {2p=0}
                         \end{Sb}
                         H_i(\pglk,\zz) \oplus
                     \bigoplus\begin{Sb}
                          {p\in \overline{E}, 2p\ne 0}\\
                           { p \sim -p}
                          \end{Sb}
                         H_i(\uk,\zz)} \\
               &   & \oplus \bigoplus\begin{Sb}
                       {l \in k, F(l,y)=0}\\ {\text{\em has no rat.
                          sol.}}
                          \end{Sb}
                         H_i(\ukw/\uk,\zz).
\end{eqnarray*}
\noindent Here, $p\sim -p$ means that we identify the factors
corresponding to the points $p$ and $-p$.}

\medskip

This is proved by examining the action of $\pgla$ on the Bruhat--Tits
tree $\X$ associated to a two-dimensional vector space over the
function field of $E$ (see, {\em e.g.,} Serre's book \cite{serre}).
The quotient graph $\pgla\backslash \X$ is a
tree (this is due to S.~Takahashi \cite{takahashi}) and the various 
simplex stabilizers are easily described.

As a corollary, we have the following rigidity result.

\medskip

{\sc Corollary.} {\em  Suppose the field $k$ is algebraically
closed and let $x,y$ be distinct points on $\overline{E}$.  Then the
specialization homomorphisms
$$s_x,s_y:H_\bullet(\pgla,\zz)\longrightarrow H_\bullet(\pglk,\zz)$$
coincide.}

\section{The Quotient and the Stabilizers}

In \cite{takahashi}, Takahashi described a fundamental domain for the
$\pgla$-action on $\X$; denote this subtree by $\D$.  There is a
distinguished vertex $o \in \D$.  For each $l$ in $k\cup \{\infty\}$ 
there is a vertex $v(l)$ adjacent to $o$.  The rest of $\D$ may be 
described as follows.  Let $\D(l)$ denote the subtree of $\D - \{o\}$ 
which contains $v(l)$.  The tree $\D$ is the union of $o$ and the
various $\D(l)$ (which are disjoint).  The trees $\D(l)$ are as 
follows.

(1) Suppose $F(x,y)=0$ has no rational solution with $x=l$.  Then 
$\D(l)$ consists only of $v(l)$ (see Figure 1).

\begin{figure}
\centerline{\psfig{figure=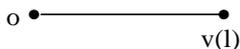,angle=270,width=1.25in}}
\caption{$F(l,y)=0$ has no rational solutions}
\end{figure}

(2) Suppose $l=\infty$ or $F(x,y)=0$ has a unique rational solution
with $x=l$.  Let $p$ be the point at infinity of $E$ or the rational
point corresponding to the solution.  Note that $p$ is a point of
order $2$.  Then $\D(l)$ consists of an infinite path $c(p,1),c(p,2),
\dots$ and an extra vertex $e(p)$ (see Figure 2).

\begin{figure}
\centerline{\psfig{figure=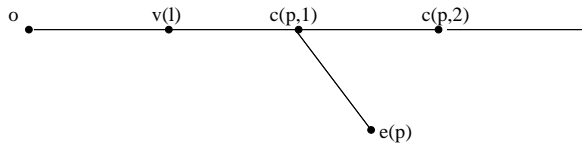,angle=270,width=3in}}
\caption{$F(l,y)=0$ has a unique rational solution}
\end{figure}

(3)  Suppose $F(x,y)=0$ has two different solutions such that $x=l$.
Let $p,q$ be the corresponding points on $E$.  Then $\D(l)$ consists 
of two infinite paths $c(p,1),c(p,2),\dots$ and
$c(q,1),c(q,2),\dots$ (see Figure 3).

\begin{figure}
\centerline{\psfig{figure=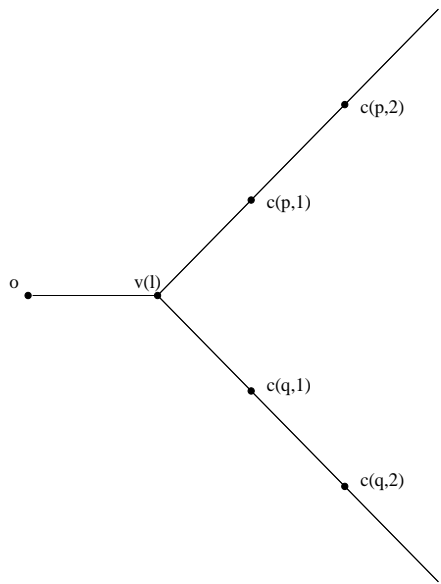,angle=270,height=3in,width=2.25in}}
\caption{$F(l,y)=0$ has two distinct solutions}
\end{figure}

The infinite path $c(p,1),c(p,2),\dots$ is called a {\em cusp}.  Note
that there is a one-to-one correspondence between cusps and the
rational points of $\overline{E}$.

Since $\X$ is contractible, we have a spectral sequence with $E^1$--term
$$E_{p,q}^1 = \bigoplus_{\sigma^{(p)} \subset \D}
H_q(\Gamma_{\sigma},\zz) \Longrightarrow H_{p+q}(\pgla,\zz)$$
where $\Gamma_{\sigma}$ is the stabilizer of the $p$--simplex $\sigma$
in $\pgla$.

The various stabilizers of the $GL_2(A)$ action were described in
\cite{takahashi}.  Denote these by $\tilde{\Gamma}_{\sigma}$.

\begin{prop}[cf. \cite{takahashi},Theorem 5]  The stabilizers
$\tilde{\Gamma}_{\sigma}$ are (up to isomorphism)
\begin{eqnarray*}
\tilde{\Gamma}_o & \cong & \uk  \\
\tilde{\Gamma}_{v(l)} & \cong & {\begin{cases}
                                 \ukw & \text{in case {\em (1)}} \\
                                 \uk \times k & \text{in case {\em (2)}} \\
                                 \uk \times \uk & \text{in case {\em (3)}}
                                 \end{cases}} \\
\tilde{\Gamma}_{e(p)} & \cong & GL_2(k) \\
\tilde{\Gamma}_{c(p,n)} & \cong & (k^n\times_{\theta} \uk)\times \uk
\end{eqnarray*}
where $k^n \times_{\theta} \uk$ is the semidirect product of $k^n$
and $\uk$ and $\theta$ is the automorphism of $k^n$ given by
coordinate-wise multiplication by elements of $\uk$.
Furthermore, the stabilizer of an edge is the intersection of the
stabilizers of its vertices.
\end{prop}

The groups $\tilde{\Gamma}_{c(p,n)}$ are of the form
$$\biggl\{\left(\begin{array}{cc}
     p  &  q  \\
     0  &  s
  \end{array}\right): p,s \in \uk, q \in k^n \biggr\}.$$
Denote the diagonal subgroups of these $\tilde{\Gamma}$ by
     $\tilde{L}$. By Theorem 1.11 of \cite{nessus}, these stabilizers
satisfy
$$H_\bullet(\tilde{\Gamma},\zz) \cong H_\bullet(\tilde{L},\zz),$$
the isomorphism being induced by the inclusion $\tilde{L} \rightarrow
\tilde{\Gamma}$.

\begin{cor} The stabilizers $\Gamma_{\sigma}$ satisfy
\begin{eqnarray*}
H_\bullet(\Gamma_o,\zz) & \cong & H_\bullet(\{1\},\zz) \\
H_\bullet(\Gamma_{v(l)},\zz) & \cong & {\begin{cases}
                                    H_\bullet(\ukw/\uk,\zz) & \text{in
case {\em (1)}} \\
                                    H_\bullet(k,\zz) & \text{in case
{\em (2)}} \\
                                    H_\bullet(\uk,\zz) & \text{in case
{\em (3)}}
                                        \end{cases}} \\
H_\bullet(\Gamma_{c(p,n)},\zz) & \cong & H_\bullet(\uk,\zz) \\
H_\bullet(\Gamma_{e(p)},\zz) & \cong & H_\bullet(\pglk,\zz).
\end{eqnarray*}
\end{cor}

\section{The Main Theorem}

Note that our spectral sequence consists of two columns and that each
row $E_{\bullet,q}^1$ is the chain complex $C_\bullet(\D,\cH_q)$,
where $\cH_q$ is the coefficient system $\sigma\mapsto
H_q(\Gamma_{\sigma})$.  Fix a positive integer $q\ge 1$.  Note that
$H_q(\Gamma_o,\zz)= 0$ and that $H_q(\Gamma_e,\zz)=0$ for any edge
incident with $o$.  It follows that the chain complex
$C_\bullet(\D,\cH_q)$ is a direct sum of chain complexes
$$C_\bullet(\D,\cH_q) = \bigoplus_{l\in k\cup \{\infty\}}
C_\bullet(\D(l),\cH_q).$$

\begin{prop}  Suppose $F(x,y)=0$ has no rational solutions with
$x=l$.  Then $$H_i(\D(l),\cH_q) = \begin{cases}
                                    H_q(\ukw/\uk,\zz) & i=0 \\
                                            0         & i>0.
                                   \end{cases}$$
\end{prop}

\begin{pf}  This is clear since $C_\bullet(\D(l),\cH_q)$ consists only
of the single group $H_q(\ukw/\uk,\zz)$ sitting in degree zero.
\end{pf}

\begin{prop}  Suppose $l=\infty$ or $F(x,y)=0$ has a unique rational
solution with $x=l$.  Then
$$H_i(\D(l),\cH_q)= \begin{cases}
                     H_q(\pglk,\zz) & i=0 \\
                                0   & i>0.
                     \end{cases}$$
\end{prop}

\begin{pf}  The stabilizer of $v(l)$ and of the edge joining $v(l)$ to
$c(p,1)$ is isomorphic to $k$.  The map $\Gamma_{\langle
v(l),c(p,1)\rangle} \rightarrow \Gamma_{v(l)}$ is an isomorphism on homology 
and the map  $\Gamma_{\langle v(l),c(p,1)\rangle} \rightarrow \Gamma_{c(p,1)}$
induces the zero map on homology.  It follows that
$H_\bullet(\D(l),\cH_q) \cong H_\bullet(\D(l)',\cH_q)$ where $\D(l)'$
is the tree obtained by deleting $v(l)$ and the edge joining it to
$c(p,1)$.  One checks easily that the map $$C_1(\D(l),\cH_q) \longrightarrow
C_0(\D(l),\cH_q)$$ is injective with cokernel $H_q(\pglk,\zz)$
(or equivalently, check that the relative homology groups 
 $H_\bullet(\D(l),e(p);\cH_q)$ vanish).
\end{pf}

\begin{prop}  Suppose $F(x,y)=0$ has two distinct solutions with
$x=l$.  Then $$H_i(\D(l),\cH_q)=\begin{cases}
                                  H_q(\uk,\zz) & i=0 \\
                                       0       & i>0.
                                \end{cases}$$
\end{prop}

\begin{pf}  In this case, $\D(l)$ is a tree and the stabilizer of each
vertex and each edge is $\uk$.  The maps $\Gamma_e \rightarrow \Gamma_v$
induce isomorphisms on homology for each edge and vertex.  It follows
that $C_\bullet(\D(l),\cH_q)$ is a chain complex with {\em constant}
coefficients.  Since $\D(l)$ is contractible, the result follows.
\end{pf}

\begin{theorem} For all $i\ge 1$,
\begin{eqnarray*}
H_i(\pgla,\zz) & \cong & {\bigoplus\begin{Sb}
                               {p \in \overline{E}}\\ {2p=0}
                             \end{Sb}
                          H_i(\pglk,\zz) \oplus
                          \bigoplus\begin{Sb}
                               {p \in \overline{E}, 2p\ne 0}\\
                               {p\sim q}
                              \end{Sb}
                          H_i(\uk,\zz)} \\
               &       & \oplus {\bigoplus\begin{Sb}
                                     {l\in k,F(l,y)=0} \\
                                     {\text{\em has no solutions}}
                                    \end{Sb}   
                                  H_i(\ukw/\uk,\zz)}.
\end{eqnarray*}
\end{theorem}

\begin{pf}  Note that since $\D$ is contractible, $E_{0,0}^2=\zz$ and
$E_{1,0}^2=0$.  The preceding propositions show that 
$$H_0(\D,\cH_q) = \bigoplus_{l\in k\cup \{\infty\}} H_0(\D(l),\cH_q)$$
and $H_1(\D,\cH_q)=0$.  It remains to identify the various direct
summands with points of $\overline{E}$.

Those $l$ for which $F(l,y)=0$ has a unique solution (or $l=\infty$)
correspond to points of order $2$ in $\overline{E}$.  Those $l$ for
which $F(l,y)=0$ has two distinct solutions correspond to pairs of
points on $\overline{E}$.  For such a pair $p,q$, the groups
$H_i(\uk,\zz)$ arising from the cusps associated to $p,q$ are
identified together since they are both adjacent to $v(l)$.  The
claimed direct sum decomposition follows.
\end{pf}

\section{Rigidity}

Suppose now that the field $k$ is algebraically closed.  Note that
the projective curve $\overline{E}$ is isomorphic to the group
$\text{Pic}^0\overline{E}$ of degree zero line bundles on
$\overline{E}$.  Moreover, after a suitable linear change of
coordinates, the points $p,q$ corresponding to the two solutions
of $F(l,y)=0$ may be assumed to satisfy $q=-p$.

\begin{cor} If $k$ is algebraically closed, then for all $i\ge 1$,
\begin{eqnarray*}
H_i(\pgla,\zz) & \cong & {\bigoplus\begin{Sb}
                              {L\in \text{\em Pic}^0\overline{E}}\\ {2L=0}
                              \end{Sb}
                              H_i(\pglk,\zz) \oplus
                          \bigoplus\begin{Sb}
                              {L\in \text{\em Pic}^0\overline{E},2L\ne 0}\\
                               {L \sim -L}
                              \end{Sb}
                              H_i(\uk,\zz)}.
\end{eqnarray*}
\end{cor}

\begin{pf}  This is obvious once one notes that the factors
$H_i(\ukw/\uk,\zz)$ do not enter the picture when $k$ is algebraically
closed.
\end{pf}

As a consequence, we have the following rigidity result.

\begin{cor}  Suppose $k$ is algebraically closed and let $x,y$ be
distinct points on $\overline{E}$.  Then the corresponding
specialization homomorphisms
$$s_x,s_y:H_i(\pgla,\zz)\longrightarrow H_i(\pglk,\zz)$$
coincide for all $i\ge 0$.
\end{cor}

\begin{pf} This follows from the direct sum decomposition of Corollary
3.1.  Since the groups which appear as summands do not involve
rational functions on $\overline{E}$, the homomorphisms $s_x$ and
$s_y$ must agree on each summand.
\end{pf}

\medskip

{\small
\centerline{\sc Acknowledgments}
I thank Andrei Suslin for pointing out that rigidity follows from
the main theorem.  I also thank Rick Jardine for many interesting
discussions about the Friedlander--Milnor Conjecture.}

\end{document}